\documentclass[11pt]{amsart}
\usepackage{graphicx,xcolor,verbatim}
\usepackage{amssymb,amsfonts,amsmath}
\usepackage{graphics,epsfig}
\usepackage{graphicx}
\usepackage{latexsym, arcs}
\usepackage{multirow}

\newcommand{\be}{\begin{equation}}
\newcommand{\ee}{\end{equation}}

\newtheorem{remark}{Remark}

\begin{document}

\title[Mathematical models of degradation for lime-based mortars]{Mathematical modelling of  the water absorption properties for historical lime-based mortars from Catania (Sicily, Italy)}

\author{G. Bretti} 
\thanks{Istituto per le Applicazioni del Calcolo ``M.Picone'', Rome, Italy, c.braun@iac.cnr.it, gabriella.bretti@cnr.it}
\author{C. M. Belfiore} 
\thanks{University of Catania, Dept. Biological, Geological and Environmental Sciences, Catania, Italy, cbelfio@unict.it}

\begin{abstract}

In this paper we propose a mathematical model of the capillary and permeability properties of lime-based mortars from the historic built heritage of Catania (Sicily, Italy) produced by using two different types of volcanic aggregate, i.e. ghiara and azolo.
In order to find a formulation for the capillary pressure and the permeability as functions of the saturation level inside the porous medium we calibrate the numerical algorithm against imbibition data. The validation of the mathematical model was done by comparing the experimental retention curve with the one obtained by the simulation algorithm. Indeed, with the proposed approach it was possible to reproduce the main features of the experimentally observed phenomenon for both materials. 
\end{abstract}


\maketitle

\section{Introduction}


 The presence of mortars in historical as well as in modern buildings is high. Their main role is to connect masonry elements such as stones and bricks, but they are also used for protective and aesthetic purposes, as  basis for floors, frescoes or mosaics. 
 
 Mortars are open porous systems that can be subject to different damaging processes due to the exposure to weathering or chemical aggression. Unfortunately, the deterioration processes of these porous materials caused by the interaction with the environment are still not completely understood, due to the complexity of phenomena involved. Moreover, their composition may vary considerably; indeed, during each historical period, clay, lime, natural pozzolans, brick dust, gypsum or combination of them were used.
 Most of deterioration processes affecting porous materials are associated to liquid water penetration, either in the form of rain or groundwater moisture infiltration, see for instance \cite{volumemach, charola} and references therein. 
 
In this respect, the assessment of the absorption properties of porous materials is crucial to determine the flow of liquid within the pores. This is not an easy task, since mineralogical composition, porosity and pore size distribution can vary even for samples of the same material.
Many mathematical models can be found in literature describing the absorption properties of porous materials, see \cite{Ass, bc, vg}.

In recent papers of one of the authors, the mathematical model describing capillary rise and water imbibition introduced in \cite{Clarelli2010} was applied with promising results to different materials, also in presence of protective treatments, see \cite{salt, goid}. 

 In the present paper a mathematical model based on Darcy's law is introduced in order to simulate numerically the water uptake into porous materials.  Since we deal with unsaturated flow, the mass balance equation for the liquid (water) 
 having density $\rho_l$ is \cite{barenblatt}:
 \begin{equation}\label{eq_richards}
\partial_t \theta = \partial_z q 
\end{equation}
where $q$ is the \emph{volumetric flux}, having the unit measure of a velocity, also called \emph{superficial velocity} of the fluid within the porous medium and it is given by the well known Darcy's law \cite{bear}: 
\begin{equation}\label{darcy}
q=-\frac{k\left(\frac{\theta_l}{n_0}\right)}{\mu_l} \left(\partial_z P_c\left(\frac{\theta_l}{n_0}\right) - \rho_l g\right).
\end{equation}
Note that $P_c$ is the capillary pressure, i.e. the pressure drop on the interface between liquid and gas, $k$ is the (intrinsic) \emph{permeability} of the porous matrix to vapour density, $\mu_l$ the \emph{viscosity} of the fluid, $n_0$ is the open \emph{porosity} of the material, i.e. the fraction of volume occupied by voids, $\theta_l$ is the \emph{fraction of volume occupied by the fluid} while the term $\rho_l g$ takes into account the effect of the gravity in the vertical flow and it can be safely disregarded for small specimens of a given material.

The equation \eqref{eq_richards} represents a continuity equation describing the variation in time of the fraction of pores occupied by the liquid due to its diffusion in the medium.
As usual in literature \cite{barenblatt, bear}, we assume the capillary pressure as a function of the fluid saturation $\theta_l/n_0$ only, i.e. $P_c = P_c(\theta_l/n_0)$.  
 
  Then, following the approach in \cite{Clarelli2010}, we introduce a function $B(s)$ depending on $s=\theta / n_0$, the \emph{water saturation}, i.e. the saturation level within the material, such that $\partial_z B(s) = -\frac{k(s)}{\mu} \partial_z P_c(s)$. 
 
 In \cite{Clarelli2010} the absorption function $B$ has a mathematical formulation of a polynomial of third degree with the derivative $B'$ resulting in a concave parabola on a compact support $s \in [s_R, s_S]$, see \eqref{NN} in paragraph \ref{sec:NN}. In this formulation, derived as an empirical/phenomenological function describing the physical properties of the porous material, there are three crucial parameters: the diffusion rate of water in the medium, the residual saturation that ensures the hydraulic continuity and the maximum value of saturation.
 
 Here, in order to describe more accurately the physical properties involved in the suction process, we introduce a new formulation of function $B$, see paragraph \ref{sec:BkP}, resulting in a sort of asymmetrical concave parabola where we define separately the capillary pressure function and the permeability function, thus allowing us to derive specific parameters related to the capillary pressure and to the permeability properties of the porous material, that can be calibrated against data. This new formulation allows us to describe by the mathematical model to curve profiles of $P_c(s)$ and $k(s)$ separately, and in this way it is possible to make a validation of the model using the available experimental data, such as  permeability test to water vapour \cite{deBoever} and/or MIP \cite{dellavecchia}. Moreover, the  model \eqref{BkP} is formulated in such a way to be more adherent to the physical process of imbibition of porous materials respect to the empirical formulation in \eqref{NN}.  
 
 Here, in order to describe the absorption properties of ghiara and azolo mortars, imbibition tests with pure water were performed in \cite{belfiore}.
 The diffusion coefficient, the residual and maximum value of saturation for the two materials are determined through the numerical calibration of the simpler model in against data - i.e. by computing the difference between model outcomes and experimental data for ghiara and azolo mortars. Then, starting from the results obtained for the simpler model, we  calibrate the six parameters of the more sophisticated model against imbibition data getting more accurate results and a complete description of the absorption properties of ghiara and azolo mortars, see Section \ref{inverse}.

The paper is organized as follows: Section \ref{sec:material} is devoted to the description of laboratory experiments on the two lime-based mortars as well as to the presentation of the mathematical model and the related numerical algorithm. Section \ref{inverse} is concerned with the definition of the inverse problem for the numerical calibration of model parameters and with the validation of the proposed model; moreover, a local sensitivity analysis on the effect of parameters on numerical results is provided. The results are discussed in Section \ref{sec:disc}. 
Finally, we end up with conclusions and future perspectives of our research in Section \ref{concl}.

\section{Materials and methods}\label{sec:material}

\subsection{The laboratory experiment \cite{belfiore}}
The laboratory experiment was carried out on two lime-based mortars replicating those of the historic built heritage of Catania (Eastern Sicily, Italy), characterized by two different volcanic aggregates, locally known as azolo and ghiara. Azolo, now obtained from the fine grinding of basalts, in its ancient meaning, was an incoherent pyroclastic rock with grain size ranging from 0.2 to 2 mm, sharp edges, and a dark-gray color \cite{belfiore2022}. Ghiara is a reddish material deriving from the thermal transformation at highly oxidizing conditions of paleo-soils, originally rich in organic matter, induced by lava flows \cite{belfiore2010, belfiore2022}.

The experimental mortar specimens have been reproduced in laboratory by following the ancient recipes known from literature \cite{battiato}. The binder used is a ready-made commercial slaked lime with CaO content higher than 90\%, MgO in the range 1.5–2\%, and solid/water ratio 40:60. The aggregate/binder ratio per volume is 2:1 in both mortars. The grain size selection of aggregates has been based on the Fuller equation by considering a maximum grain size of 2 mm (mortars for finishing purposes).
The two replicated mortars were designed in \cite{belfiore} to assess their physical-chemical durability. For that purpose, the authors used a multi-analytical approach including mineralogical-petrographic and physical investigations (polarized optical microscopy, X-ray diffraction, pore size distribution, water absorption by capillarity, water vapor permeability), along with accelerated aging tests by salt crystallization and sulfur dioxide. The data of capillary water absorption test, here used for the mathematical modelling, were determined according to the European standard UNI EN 1015-18 (2004). Specifically, three specimens per each mortar type were placed inside a ventilated oven for 24 h to remove residual moisture and after that the dry weight of each sample was calculated. Then, the specimens were placed inside a container with the lower face immersed in about 5–10 mm of water and covered to prevent water evaporation. After 10 min, the specimens were extracted from the basin and weighted (M1), then again immersed and after 90 min weighted to obtain the final mass (M2). The capillary coefficient was calculated through the following equation:

$$C= 0.1(M_2 - M_1) [g/ (cm^2 s^{1/2})].$$

\begin{figure}
\centering
\includegraphics[height=7 cm, width=8.5cm]{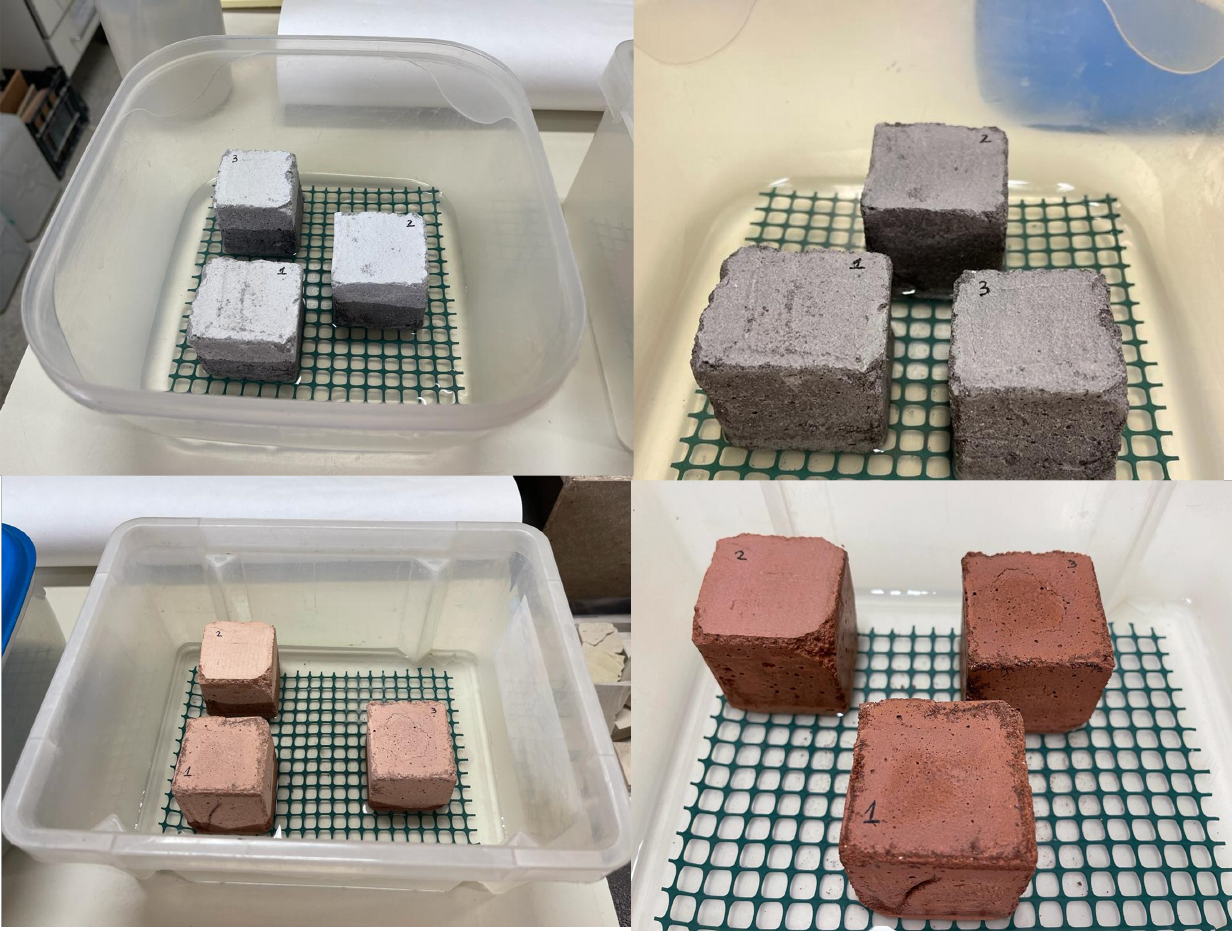}
\caption{Laboratory experiment for the measurement of water absorption for azolo and ghiara mortars.}  
\label{fig:assorb}   
\end{figure}

The amount of absorbed water ($Q_i$) was calculated as follows:

$$Q_i = (w_i - w_0)/A,$$

where: $w_i$ and $w_0$ are the weights (in g) of the sample at times $t_i$ and $t_0$, respectively; A ($cm^2$) is the area of the surface exposed to water.
See Fig. \ref{fig:assorb} for the experimental setup of the imbibition test.

\subsection{MIP experiment}\label{sec:MIP}
The experimental test used to have quantitative information on  the pore size distribution of porous media with interconnected porosity is the mercury intrusion porosimetry (MIP). Specifically, pores in the range 0.005–750 $\mu m$ were analyzed by  using a Micromeritics Autopore V 9600 porosimeter reaching a maximum pressure of 227 MPa \cite{belfiore}. 

By using information provided by MIP, water retention models can be developed, see \cite{dellavecchia} and references therein.

 Suction is computed from MIP data using Laplace formula reported in \cite{CVB}. In particular, the capillary pressure, $P^i_{MIP}$, can be computed directly from mercury pressure, $P^i_{Hg}$ corresponding to saturation level $s^i_{Hg}$, using the Laplace equation:
\begin{equation}\label{laplace}
    P^i_{MIP} = \frac{T_w \cos(\theta_w)}{T_{Hg} |\cos(\theta_{Hg})|} P^i_{Hg},
\end{equation}
with $i=1,\ldots, N_{MIP}$ the number of MIP measurements, $T_w$ and $T_{Hg}$ the surface tensions of water and mercury, respectively ($T_w = 0.073 N/m$ and $T_{Hg} = 0.489 N/m$); and $\theta_w$ and $\theta_{Hg}$ the contact angles of water and mercury, respectively, $\theta_w=0^\circ$ and $\theta_{Hg}=130^\circ$. We remark that here we express the capillary pressure in the unit measure $[g/(cm \cdot s^2)].$

Notice that the saturation level $s_{Hg,i}$ is obtained from MIP data as:
$$
s_{Hg,i} = 1-V_i/V_{max}, \ i=1,\ldots, N_{MIP},
$$
with $V_i$ and $V_{max}$, respectively, the specific and the maximum specific volume of pores.

\subsection{The mathematical model of suction of pure water}
We denote the fraction of volume occupied by the liquid and by the gas (composing the fluid) within the representative element of volume, respectively by $\theta_l$ and $\theta_g$. The open porosity $n_0$ of unperturbed material satisfies the following relation:
\begin{equation}\label{porosity}
n_0 = \theta_l + \theta_g.
\end{equation}
From now on, for the sake of simplicity, we suppress the label $l$ and, since we are considering specimens of small dimension, gravity effects can be safely neglected in Darcy’s law. Then, we introduce a function $B$ describing the absorption properties of porous material so that Darcy's law \eqref{darcy} can be rewritten as 
\begin{equation}\label{Bdarcy}
q=\partial_z B(s) = -\frac{k(s)}{\mu} \partial_z P_c(s).
\end{equation} 

The identification of the function $B$ reproducing the capillary rise properties of a specific porous material is not an easy task.
Indeed, there are many suggested experimental curves aiming to connect the capillary pressure with the moisture content, but they are not completely reliable since a relation correlating capillary pressure with moisture content
into the porous matrix still lacks. There are many possible choices for the
function $B(s)$, since it can vary very much according to the physical property
of any given material.
A simple and phenomenological description of the capillary rise of water can be obtained by defining $B$ as a compactly supported function  in order to reproduce the fact that the liquid (water) flows within the material only if the saturation $s$ ranges between a minimum and maximum saturation level, with a diffusion rate to be calibrated against experimental data.  Two possible formulations of $B'$, i.e. "concave parabola" and "asymmetrical concave parabola" are proposed, respectively, in the next paragraphs \ref{sec:NN} and \ref{sec:BkP}, see also pictures in Fig. \ref{fig:BNN}.

 Then, the one-dimensional mathematical model for the experiment of imbibition given by \eqref{eq_richards} can be written as: 
\begin{equation}\label{pb-water}
\partial_t \theta = \partial_{zz} B(s),
\end{equation}
with $B$ the function satisfying relation \eqref{Bdarcy}.  

Equation (\ref{pb-water}) has to be coupled with reasonable initial and boundary conditions.
We assume the following initial and boundary conditions 
\begin{align}
\theta(z,0)&=0, z \in [0, h_1]\label{ic}\\
\theta(x,t)&=n_0, x \in [0,h_2], t \in [0, T_f]\label{bc1}\\
\theta(h_1,t) &= \theta_{ext}, t \in [0, T_f] \label{bc2}, 
\end{align}
meaning that the bottom side of the sample is always saturated since it is placed in direct contact with water, while at the top side $z=h_1$ the condition \eqref{bc2} reproduces the exchange of the specimen with the humidity within the bucket ($\theta_{ext}$ is the moisture content of the ambient air) due to evaporation, thus we assume it as the value of humidity on the facelet on the top of the material.
In the above condition, $\theta_{ext}$ is assumed constant during the experiment.

\begin{remark}
    {The Dirichlet boundary condition at the top boundary can be replaced by the Robin boundary condition:
$$
\partial_z \theta(L,t) = K_w (\theta(L,t) - \theta_{ext}),
$$
in order to reproduce a loss of water occurring gradually at the upper boundary $z=L$, being regulated by 
the condition on $\partial_z \theta$, with a rate given by a coefficient $K_w>0$. Even if this condition seems to be more realistic, here we prefer to use a Dirichlet boundary condition in order to reduce the number of model parameters to be calibrated. However, for large values of $K_w$ the solutions obtained numerically assuming Robin boundary conditions seem to converge towards those with Dirichlet ones.  }
\end{remark}

Note that in this experimental settings the evaporation from lateral sides is negligible, thus in this in one-dimensional formulation we did not include it in \eqref{Bdarcy}.
Note that to get $\theta_{ext}$ we used formula for the saturated vapor density (SVD in [$g/m^3$ ]) as function of temperature T [ C]: 
$$ SVD(T) = (5.018 + 0.32321 T + 8.1847 \times 10^{-3} T^2 + 3.1243 \times 10^{-4} T^3) \times 10^{-6},$$
 from which we can obtain the density of vapor in [$g/cm^3$ ].

The value of $\theta_{ext}$ is obtained from the formula \cite{hyperphy}:
\begin{equation}\label{thetabar}
\theta_{ext} = \frac{SVD(T)}{\rho} \cdot n_0\cdot UR,
\end{equation}
where $T=25^{\circ}C$ and $UR=90\%$ is the percentage of humidity in the ambient air, see Table \ref{table:1} for the values of the parameters of the problem.

In the following, we propose two possible formulations for B function.

\subsection{The absorption function $B'$ \cite{Clarelli2010}}\label{sec:NN}
\begin{figure}
    \centering
\includegraphics[width=0.4\linewidth]{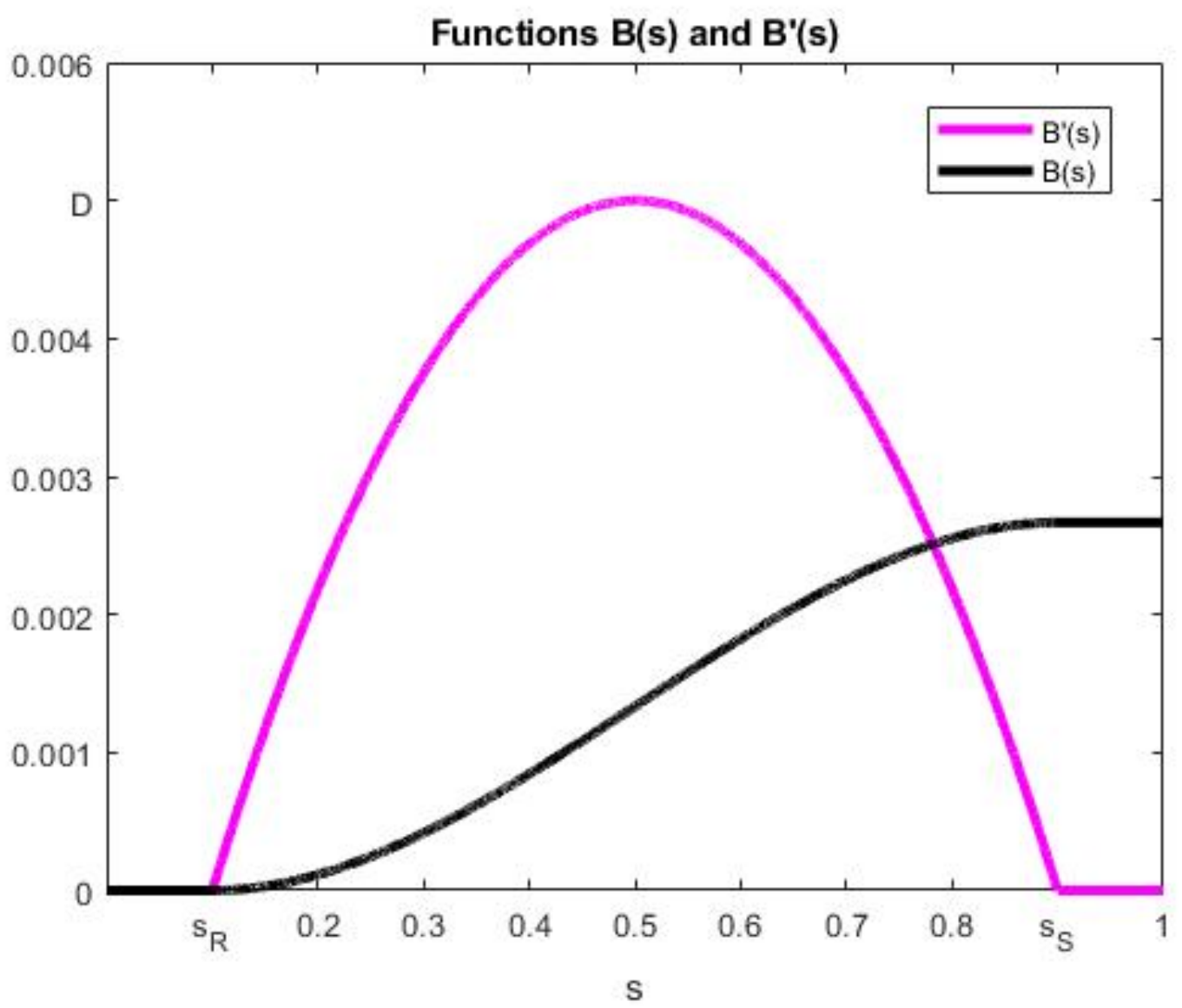}\quad \includegraphics[width=0.4\linewidth]{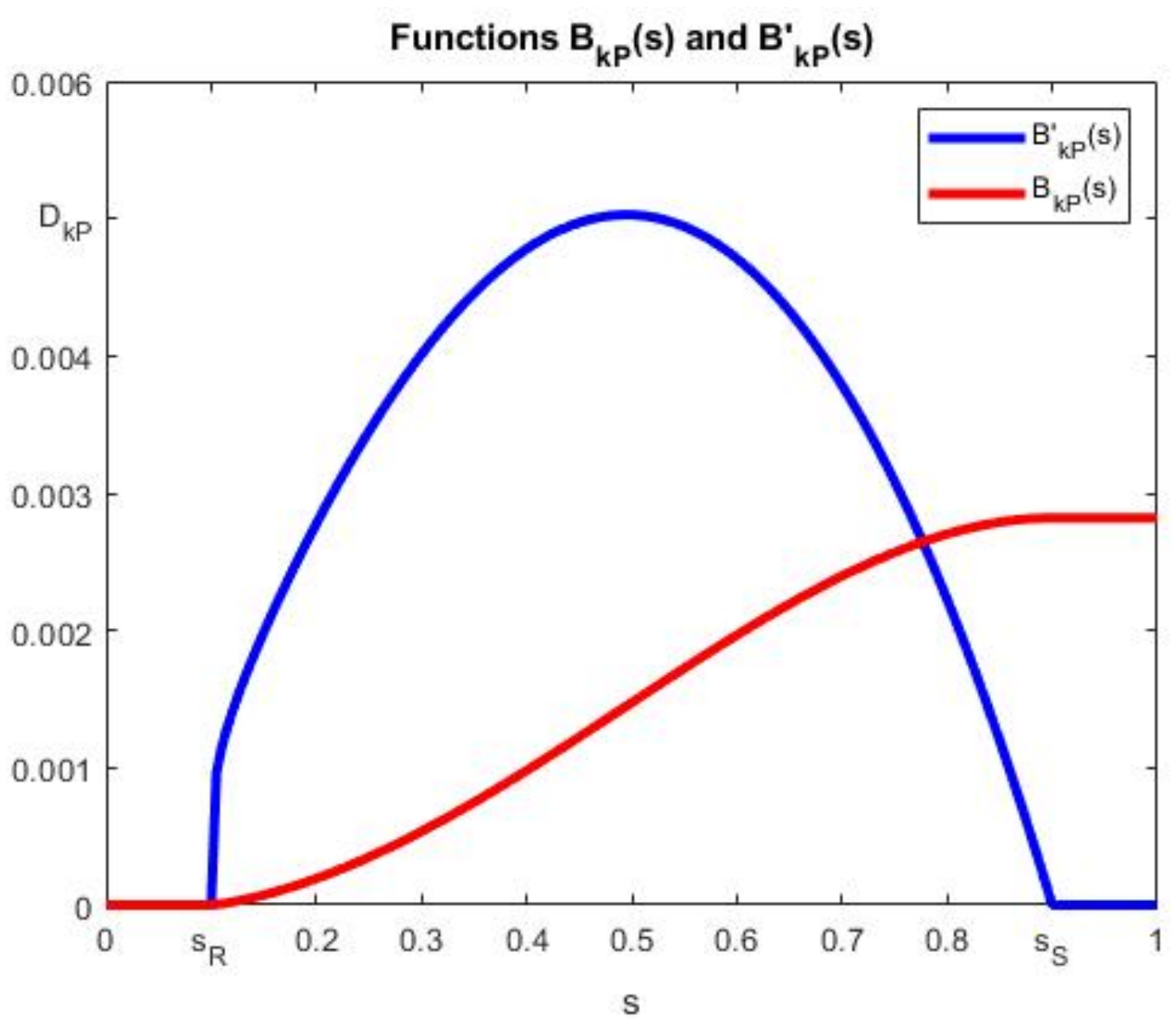}
    \caption{Left panel: Plot of functions $B(s)$ and $B'(s)$ defined in \eqref{NN} and obtained for $s_R=0.1, s_S=0.9$ and $D=$5e-03. Right Panel: Plot of functions $B(s)$ and $B'(s)$ defined in \eqref{BkP} and obtained for $s_R=0.1, s_S=0.9, \alpha=0.25, c=$1.35e+05, $K_s=$ 8e-10,$\gamma=1.45$ with $D_{kP}=$5e-03. }
    \label{fig:BNN}
\end{figure}

 Here, among the possible choices for $B$ function, a possible formulation is the function first introduced in \cite{Clarelli2010}:
\begin{align} \label{NN}
    B(s) &= 
    \begin{cases}
    0 & s \in [0,s_R)  \\  
    -\frac{(2 D (s_R - s)^2 (s_R - 3 s_S + 2 s))}{3 (s_R - s_S)^2} & s \in [s_R,s_S] \\
     \frac{2}{3} D (s_S - s_R) & s > s_S,
     \end{cases} \\
    \textrm{ with } \partial_s B(s) = B'(s) &= max\left(0,-\frac{4 D (s_R - s) (s_S - s)}{(s_R - s_S)^2}\right),
\end{align}

and $\{s_R, s_S, D\}$ is the set of model parameters to be determined: $s_R$, the minimum value for saturation ensuring the hydraulic continuity, $s_S$ the maximum value of $s$ reached at saturation, and the maximum value $D=\max_{[s_R,s_S]}{B'}$ reached at $s=\frac{s_R+s_S}{2}$. We remark that $B'(s)$ is a compactly supported function for $s$ in $[s_R,s_S]$ and $D$ has the dimensions of a diffusivity. 
Then, the profiles of $B(s)$ and $B'(s)$ are obtained by letting the right and left endpoints $s_R$ and $s_S$ vary, i.e. $s_S>s_R>0, s_S\le 1$, meaning that the realistic pore saturation is normally less than $100\%$. An example of the profiles of $B$ and $B'$ can be found in the left panel of Fig. \ref{fig:BNN}. 

\begin{remark}It is worth noting that, due to its simplicity, this formulation of $B$ and $B'$ has the great advantage of depending only on 3 parameters, thus it is quite easy to be calibrated against data. However, here we do not have the possibility of expressing the capillary pressure $P_c$ and the permeability function $k(s)$ independently as in Darcy's law \eqref{Bdarcy}. From a modelling point of view, the model \eqref{NN} can be considered as a phenomenological/empirical model. 
\end{remark}

\subsection{A new formulation of absorption function $B'$ for Darcy's law}\label{sec:BkP}
A slightly different formulation of function $B$ in Darcy's law \eqref{Bdarcy} that conserves the main features of function $B$ described in \eqref{NN} and allows us to express separately the permeability function $k(s)$ and the capillary pressure $P_c(s)$ is introduced. In particular, the permeability function $k(s)$ is formulated in such a way to reproduce the profile described in Chapter 9 of \cite{bear}, i.e the permeability is a monotonic function increasing with $s$. Then it is expressed as:
\begin{equation}\label{perm_fun}
 k(s) = \left\{
\begin{array}{ll}
K_{s} \left(\frac{s-s_R}{s_S-s_R}\right)^\gamma, \textrm{ if } s \in (s_R, s_S),\\
0 ,\textrm{ if } s \in [0,s_R],\\
K_s, \textrm{ if } s \in [s_S,1],  
\end{array}\right.
 \end{equation}
 with $K_{s}>0$ the constant of permeability at saturation, possibly obtained by experimental measurements, and $\gamma>0$ a parameter to be calibrated with experimental data. 
For the capillary pressure, a possible formulation qualitatively reproducing the profile reported in \cite{bear, dellavecchia}, 
is:
\begin{equation}\label{Pc}
P_c (s) =  c \frac{(s-s_S)^2}{(s-s_R)^{\alpha}}, \ \textrm{ if } s \in (s_R,s_S],
\end{equation}
with $c>0$ a characteristic coefficient for the given material, $s_R$ and $s_S$ the saturation parameters defined above and $\alpha \in (0,1)$ an exponent to be calibrated against data.  More in details, function $P_c(s)$ has the following behavior:  $P_c(s_S)=0$ and it tends to infinity for $s \to s_R$. 
The derivative $\partial_s P_c=P'_c (s)$ is then:
\begin{equation}\label{P'}
P'_c (s) = -\frac{c (s-s_S) (2 s_R - 2s - \alpha s_S + \alpha s)}{(s - s_R)^{\alpha+1}}.
\end{equation}

\begin{figure}
    \centering
    \includegraphics[width=0.4\linewidth]{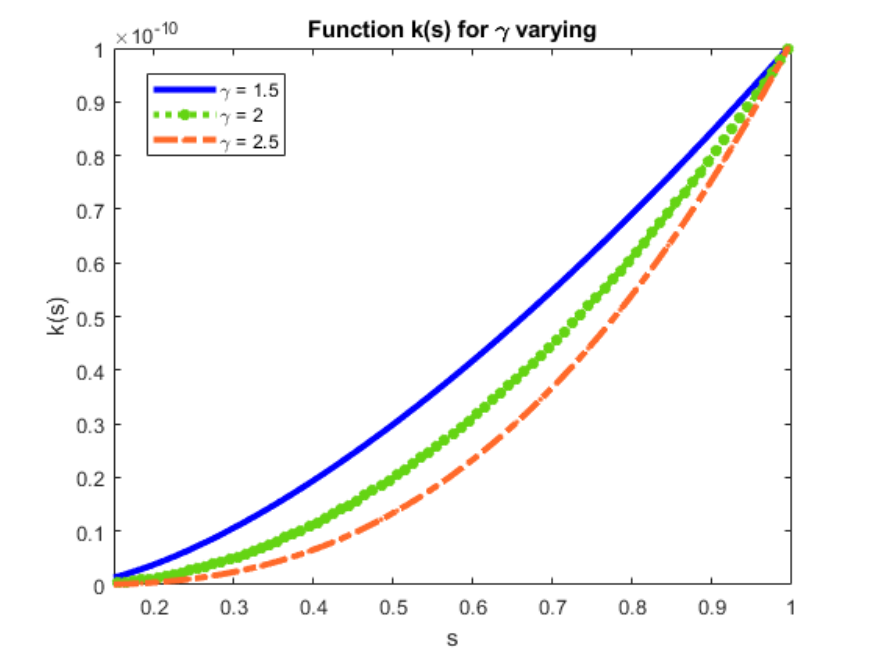}\quad  \includegraphics[width=0.4\linewidth]{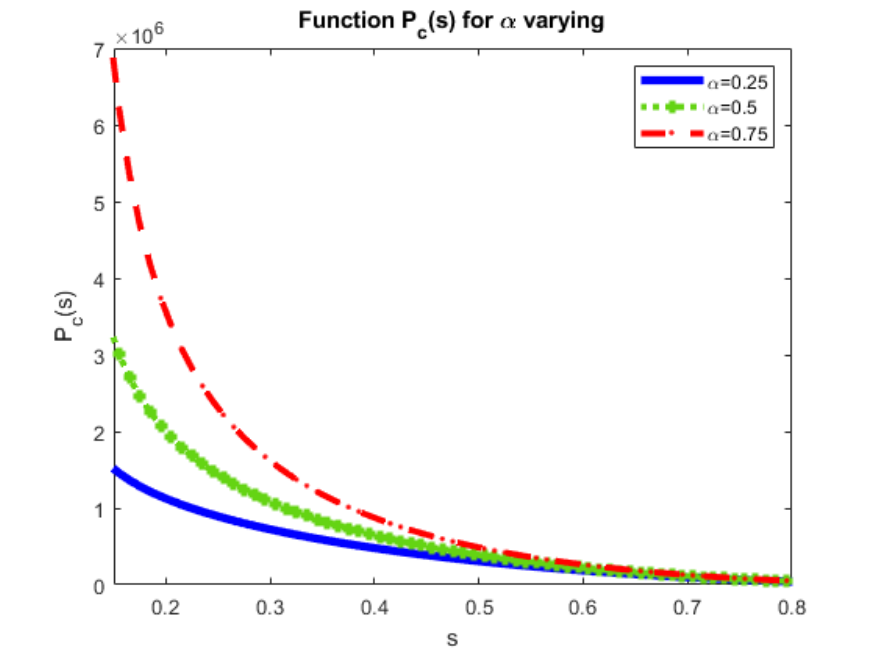}
    \caption{Plot of permeability function \eqref{perm_fun} for $K_s=$1e-10 and of capillary function \eqref{Pc} for fixed $ c=$1e+06, assuming different values of $\alpha$.}
    \label{fig:k_P}
\end{figure}

As an example, we depict in Fig. \ref{fig:k_P} the profile of functions \eqref{perm_fun} and \eqref{Pc} obtained for different parameters.

Then, keeping the same framework introduced in the previous paragraph we consider $B'$ as a compactly supported function in $[s_R, s_S]$. We introduce this new formulation using \eqref{perm_fun} and \eqref{P'} in \eqref{Bdarcy} so that we get the following expression for the derivative $\partial_s B_{kP}$ called for simplicity $B'_{kP}(s)$:
\begin{align} \label{Bpgen}
 B'_{kP}(s) &= 
    \begin{cases}
    0 & s \in [0,s_R] \cup [s_S,1]\\  
    max\left(0, K_s \frac{c}{\mu} \frac{(s-s_R)^{\gamma - \alpha -1}}{(s_S-s_R)^\gamma} (s-s_S) (2 s_R + s (\alpha -2) - \alpha s_S) \right) & s \in (s_R,s_S) .
     \end{cases} 
     \end{align}

 In this case the set of model parameters is given by $\{s_R, s_S, \alpha,c,K_{s},\gamma\}$, with $\gamma-\alpha-1>0$ and we indicate by $D_{kP} = \max_{[s_R,s_S]}{B'_{kP}}$ the diffusion coefficient. Now, we integrate exactly $B'_{kP}(s)$ from \eqref{Bpgen}, and we get the following expression:
 \begin{align} \label{BP2}
     \widetilde{B}_{kP} (s) &= 
     \frac{K_s c (s - s_R)^{\gamma - \alpha}}{\mu (s_S - s_R)^{\gamma}} \cdot\\
     &\cdot \frac{s^2 u + s v + \gamma^2 s_S (-2s_R + \alpha s_S) + \gamma z + \alpha s_S^2(\alpha^2 - 3 \alpha + 2)}{ (- \alpha^3 + 3 \alpha^2 (\gamma+1) - 3\alpha \gamma(\gamma+2) - 2\alpha + \gamma^3 + 3 \gamma^2 + 2 \gamma)},  \nonumber
     \end{align}
     with 
     \begin{align*} 
     u&= \alpha^3 - \alpha^2(2\gamma +3) +\alpha(5 \gamma +2) - 2\gamma(\gamma+1),\\
     v&=  2\gamma(- s_R \alpha + 2 \alpha^2 s_S + 2 s_S - 4\alpha s_S) +  2\gamma^2 (s_R - \alpha s_S + s_S ) + 2 \alpha s_S (-\alpha^2 + 3 \alpha - 2) ,\\
     z&= 2 s_R^2 + \alpha s_S^2 (-2\alpha + 3) + 2s_R s_S(\alpha -2),
     \end{align*}
     and we introduce the coefficient $C= K_s c$.
 Then, $B_{kP}(s)$ reads as:
\begin{align} \label{BkP}
    B_{kP}(s) &= 
    \begin{cases}
    0 & s \in [0,s_R] \\  
    \widetilde{B}_{kP}(s) \textrm { in \eqref{BP2}} & s \in (s_R,s_S) \\
      \widetilde{B}_{kP}(s_S)=\frac{2 K_s c \gamma (s_S - s_R)^{2-\alpha}}{\mu (- \alpha^3 + 3 \alpha^2(\gamma + 1) - 3 \alpha \gamma(\gamma +2) - 2\alpha + \gamma^3 + 3\gamma^2 + 2\gamma)} & s \in [s_S,1].
     \end{cases} 
     \end{align}


An example of the profiles of $B_{kP}(s)$ and $B'_{kP}(s)$ can be found in the right panel of Fig. \ref{fig:BNN}. 
\section{Calibration and validation of the mathematical model: statement of the inverse problem and results}\label{inverse}

Now we describe the calibration procedure to determine parameters of the imbibition function \eqref{Bpgen} using the experimental data of imbibition experiment both for ghiara and azolo mortars.
The objective is to reproduce numerically the cumulative water absorption, i.e. the water that crosses the mortar samples
by capillary suction.
\subsection{Calibration of imbibition parameters in \eqref{BkP} against data }
\label{sec:calib}
We calibrate model parameters in two steps:
\begin{itemize}
    \item[Step 1.] First we calibrate model parameters, i.e. saturation parameters $s_R, s_S$ and the diffusion coefficient $D$, in the simpler model \eqref{NN} against data by finding those that minimize the error \eqref{err_fun}. In this way we gain insights on the order of magnitude of the diffusive physical phenomenon corresponding to $\max_{[s_R,s_S]}{B'(s)}$ and the initial guess for $s_R, s_S$.
    \item[Step 2.] {Then, for the calibration of the model \eqref{BkP} we look for the set of parameters minimizing the difference \eqref{err_fun} between data and simulation outcomes proceeding as follows:
    \begin{itemize}
        \item we assume the values obtained at the end of Step 1 for $s_R$ and $s_S$ as initial guess and we let them vary slightly;
        \item we tune the other model parameters in $B_{kP}$, namely $\alpha, c, K_s$ and $\gamma$, in such a way that $D_{kP}$ has the same order of magnitude of the diffusion coefficient $D$ obtained in Step 1.
    \end{itemize} 
   }
    \end{itemize}
    
    Finally, we validate our results by comparing the capillary pressure profile in \eqref{Pc} using optimal parameters found with the procedure above with MIP data described in \ref{sec:MIP}.

\vspace{0.5cm}
\underline{Definition of the error between data and model outcome.} {In order to calibrate parameters of adsorption function $B'$ we need a measure to establish a comparison between model outcomes and data coming from the experiment of water imbibition. Since the relevant quantity that is experimentally accessible is the total quantity of water $Q_k$ absorbed by the specimen at certain
scheduled time intervals $t_k$, we compute the total amount of fluid absorbed obtained by the numerical algorithm at the same time intervals  ($Q^{num}_k$), as explained below.}
{First we make a discretization in time and space of the computational domain: $z_j= j\Delta z, j=0,...,N=\left[\frac{h_1}{\Delta z}\right] $, $\{t_k\}_{k=1,...,N_{meas}}$, and we define the water content on the numerical grid as $\theta^k_j =\theta(z_j, t_k)$. Then, we  solve numerically equation (\ref{pb-water}) coupled with the boundary conditions \eqref{bc1}-\eqref{bc2} in order to get the water content, by using the explicit forward-central approximation scheme \cite{diele}:
$$
\theta^{k+1}_j = \theta^k_j + \frac{\Delta t}{{\Delta z}^2} (B(\theta^k_{j+1}/n_0)-2B(\theta^k_j/n_0)+B(\theta^k_{j-1}/n_0)),
$$
where $B$ is the absorption function either defined in \eqref{NN} or \eqref{BkP}} under the CFL condition 
$$
\frac{\Delta t}{{\Delta z}^2} \le \frac{n_0}{2 \max_{[s_R,s_S]}{B'(s)}}.
$$
 At the bottom boundary we assume the imbibition condition
$\theta^{k+1}_0 = n^{k+1}_0,$
and at the top boundary we use the exchange condition: $\theta^{k+1}_{N+1} = \bar\theta$.

Then, we get the quantity of water in the specimen at time $t_k$ obtained by the numerical algorithm by approximating the integral 
$\int_{0}^{h_1} \rho \theta(z,t_k) dz,$
with the trapezoidal rule:
$$
{Q^{num}_k} = \rho\frac{\Delta z}{2} \left(\theta^k_0 + 2\sum_{j=1}^{N-1} \theta^k_j + \theta^k_{N} \right),
$$
in order to compare model outcome with experimental data at the same times. 
The error to be minimized is then defined as the sum of all the relative square errors between experimental data and numerical simulation:
\begin{equation}\label{err_fun}
E_{1,2} = \frac{1}{N_{meas}}\sum_{k=1}^{N_{meas}} \frac{
(Q^{num}_k - Q_k)^2}{(Q^{num}_k)^2},
\end{equation}
 with $E_1$ depending on the set of parameters $\{s_R, s_S, D\}$ and $E_2$ depending on $\{s_R, s_S, \alpha, c, K_s,\gamma\}$, for the functions $B'$ and $B'_{kP}$, respectively.
{The calibration procedure has been carried out in MATLAB\textcircled{c} applying the simulated annealing method; the optimal parameters obtained  for the absorption functions \eqref{NN} and \eqref{BkP} are reported in Table \ref{table:2} and \ref{table:3}, respectively}. The computational time for a single simulation with fixed parameters takes about 3 seconds for the simpler model \eqref{NN} and 120 seconds for the model \eqref{BkP} on an Intel(R) Core(TM) i7-3630QM CPU 2.4 GHz assuming as space step $\Delta z =$ 2.5e-02 and time step $\Delta t= 3.5$e-02. In Table \ref{table:err} are listed the results obtained by the absorption function \eqref{NN} and \eqref{BkP}.

\begin{table*}[tbp]
\centering
{
\begin{tabular}{|p{2cm}|c|c|c|c|} \hline
 Parameter  &  Description& Units & Value& Ref.\\\hline\hline
$n_0$ & open porosity of ghiara mortar & - & 4.66e-01 & exp. data\\
$n_0$ & open porosity of azolo mortar & - & 3.85e-01 & exp. data\\
$\tau$ & tortuosity factor of ghiara mortar & - & 9.9 & exp. data\\
$\tau$ & tortuosity factor of azolo mortar & - & 7.6 & exp. data\\
\hline\hline
\end{tabular} }
\vspace{0.2 in}   
\caption{Physical parameters of the porous materials.}\label{table:phys}
\end{table*}

\begin{table*}[tbp]
\centering
{
\begin{tabular}{|p{2cm}|c|c|c|c|} \hline
 Parameter  &  Description& Units & Value& Ref.\\\hline\hline
$\bar\theta_l$ & moisture content of the ambient air & - & 2.1e-03 & eq. \ref{thetabar} \\
$h_1$  & sample's height& cm &  5.0 & exp. data\\
$h_2$ & sample's height immersed & cm &2.5e-02& assumption\\
$\rho_l$ & density of water & $g/cm^3$ & 1  &\cite{densacqua}\\
$\mu$ & water viscosity at 25$^\circ$ & Poise&
8.9e-03& \cite{visc}\\
$T_f$ & final observation time & s & 5400& exp. data\\
\hline\hline
\end{tabular} }
\vspace{0.2 in}   
\caption{Fixed parameters/experimental setup given in input of the mathematical model \eqref{pb-water}-\eqref{bc1}-\eqref{bc2}.}\label{table:1}
\end{table*}



\vspace{1 in} 

\begin{table*}[tbp]
\centering
{
\begin{tabular}{|p{2cm}|c|c|c|} \hline
\hline
 \multicolumn{4}{|c|}{Ghiara mortar} \\
 \hline
 Parameter  &  Description& Units & Value\\\hline\hline
$s_R$ & residual saturation  &  - &  6.75e-01 \\
$s_S$ & max. saturation level  &  - &  1e-01  \\
$D$ & diffus. coeff. in \eqref{NN}& $cm^2/s$ &  1.95e-02 \\\hline\hline
\multicolumn{4}{|c|}{Azolo mortar} \\
 \hline
 Parameter  &  Description& Units & Value\\\hline\hline
$s_R$ & residual saturation  &  - &  5.21e-01 \\
$s_S$ & max. saturation level  &  - &  1e-01 \\
$D$ & diffus. coeff. in \eqref{NN}&  $cm^2$/s & 4.8e-03 \\
\hline\hline
\end{tabular} }
\vspace{0.2 in}    
\caption{Parameters obtained with a fitting procedure for $B$ function \eqref{NN}.}\label{table:2}
\end{table*}

\begin{table*}[tbp]
\centering
{
\begin{tabular}{|p{2cm}|c|c|c|} \hline
\hline
 \multicolumn{4}{|c|}{Ghiara mortar} \\
 \hline
 Parameter  &  Description& Units & Value\\\hline\hline
$s_R$ & residual saturation  &  - &  6.75e-01 \\
$s_S$ & max. saturation level  &  - &  9.994e-01 \\
$c$ & characteristic coeff. & $g \ cm^{-1} s^{-2}$  &  1.4e+06  \\
$\alpha$ & exponent in \eqref{BkP}  & -& 0.25\\
$K_s$ & permeability at satur.&  $cm^2$ & 7.65e-10\\
$\gamma$ & curvature parameter in \eqref{BkP} &  - & 1.865 \\
$D_{kP}$ & diff. coeff. in \eqref{BkP} &  $cm^2 {s}^{-1}$&  1.97e-02\\
\hline\hline
 \multicolumn{4}{|c|}{Azolo mortar} \\
 \hline
 Parameter  &  Description& Units & Value\\\hline\hline
$s_R$ & residual saturation  &  - &  5.5e-01 \\
$s_S$ & max. saturation level  &  - &  1e-01 \\
$c$ & characteristic coeff. & $g \ cm^{-1} s^{-2}$ &  1.98e+05\\
$\alpha$ & exponent in \eqref{BkP} & -& 0.25\\
$K_s$ & permeability at satur. &  $cm^2$ & 7.93e-10 \\
$\gamma$ & curvature parameter in \eqref{BkP} &  - & 1.45 \\
$D_{kP}$ & diff. coeff. in \eqref{BkP} &  $cm^2 s^{-1}$ &  4.84e-03\\
\hline\hline
\end{tabular} }
\vspace{0.2 in}     
\caption{Parameters obtained with fitting procedure for $B_{kP}$ function \eqref{BkP}.}\label{table:3}
\end{table*}

\begin{table*}[tbp]
\centering
{
\begin{tabular}{|p{3cm}|c|c|} \hline
Material  &  Error Value $E_1$ &  Error Value $E_2$\\ \hline\hline
ghiara mortar & 2.20e-04 & 1.45e-04 \\
azolo mortar &  1.44e-04 & 1.22e-04\\ 
\hline\hline
\end{tabular} }
\vspace{0.2 in}    
\caption{Error values $E_1$ and $E_2$ defined in \eqref{err_fun}  and obtained, respectively, by the absorption function \eqref{NN} and \eqref{BkP}.}\label{table:err}
\end{table*}
\begin{figure}
    \centering
    \includegraphics[scale=0.4]{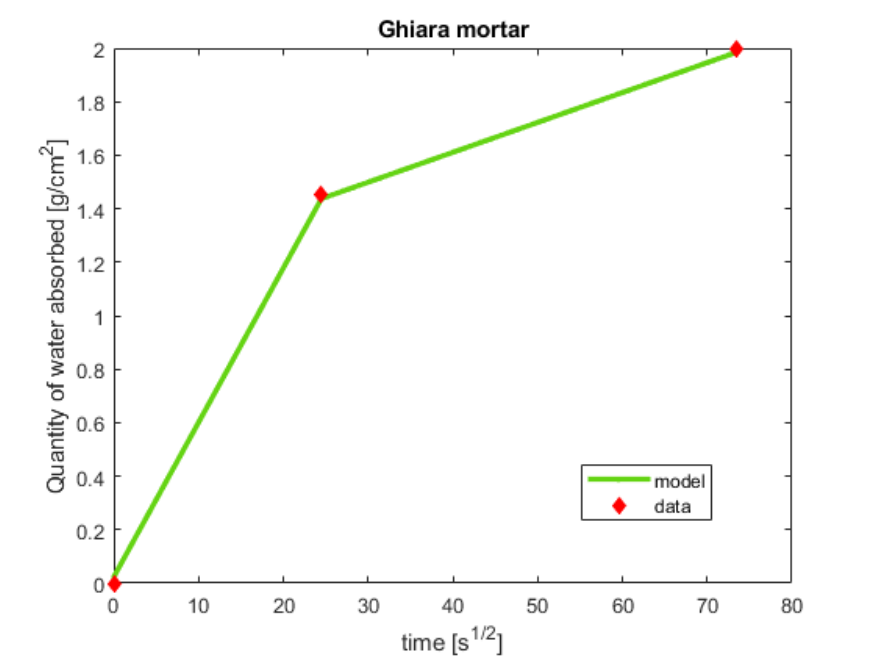}\quad \includegraphics[scale=0.4]{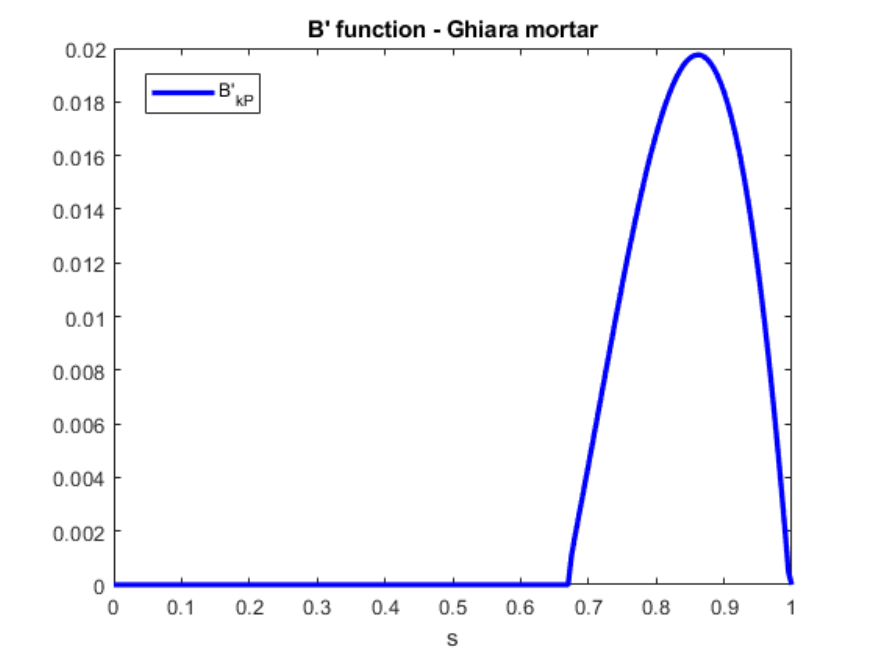}
    \caption{Left panel: Plots of the imbibition curve for ghiara mortar obtained numerically by the mathematical model \eqref{pb-water}-\eqref{bc1}-\eqref{bc2}-\eqref{BkP} (green line) VS experimental data (red diamonds). Right panel: Plot of $B'$ function (\ref{BkP}) for ghiara mortar (blue line).}
    \label{fig:ghiara}
\end{figure}
\begin{figure}
    \centering
    \includegraphics[scale=0.4]{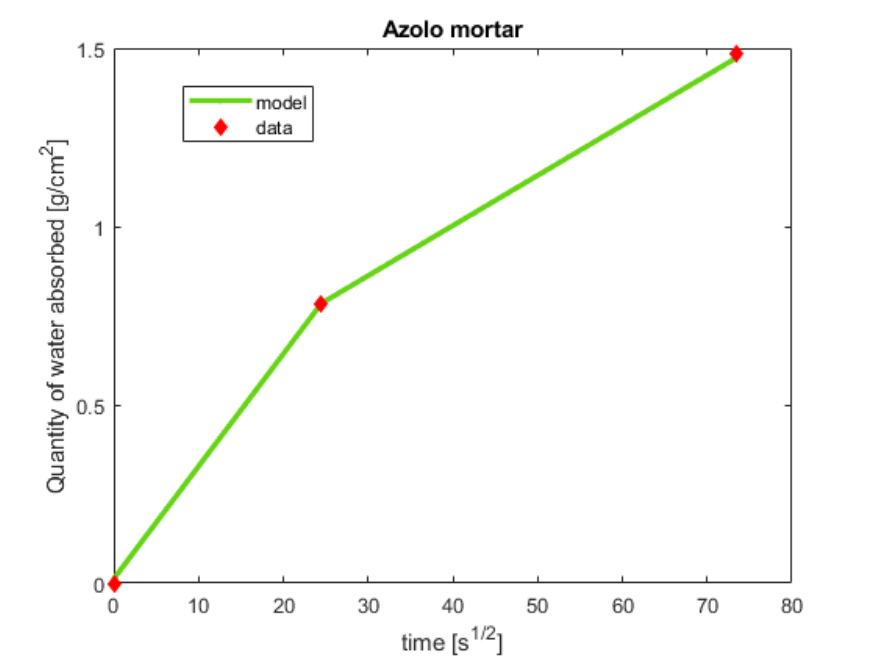}\quad \includegraphics[scale=0.4]{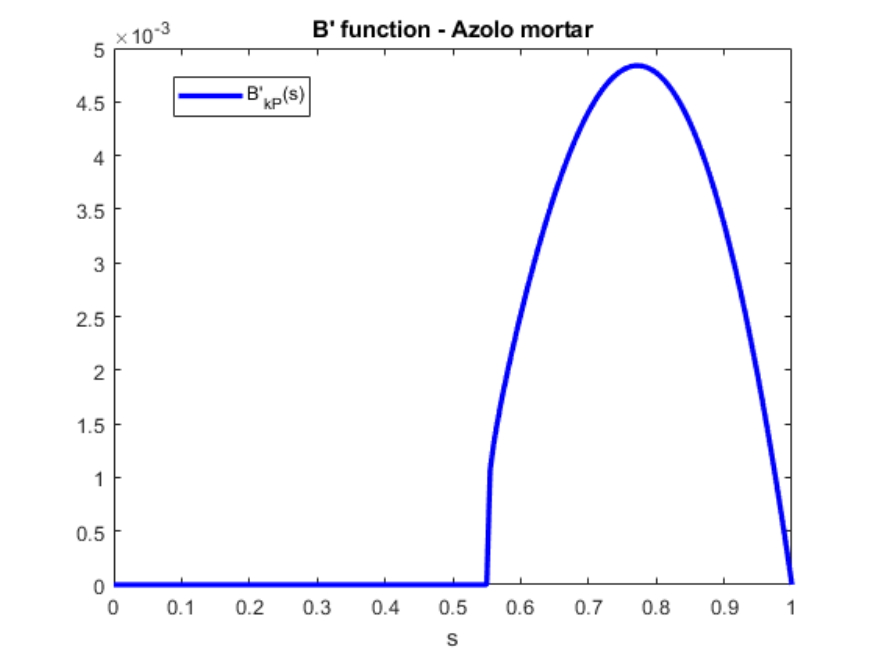}
    \caption{Left panel: Plot of the imbibition curve for azolo mortar obtained numerically by the mathematical model \eqref{pb-water}-\eqref{bc1}-\eqref{bc2}-\eqref{BkP} (green line) VS experimental data (red diamonds). Right panel: Plot of $B'$ function  (\ref{BkP}) for azolo mortar (blue line).}
    \label{fig:azolo}
\end{figure}

\subsection{Validation of the calibration procedure using MIP data}\label{sec:valid}

\begin{figure}
    \centering
    \includegraphics[scale=0.4]{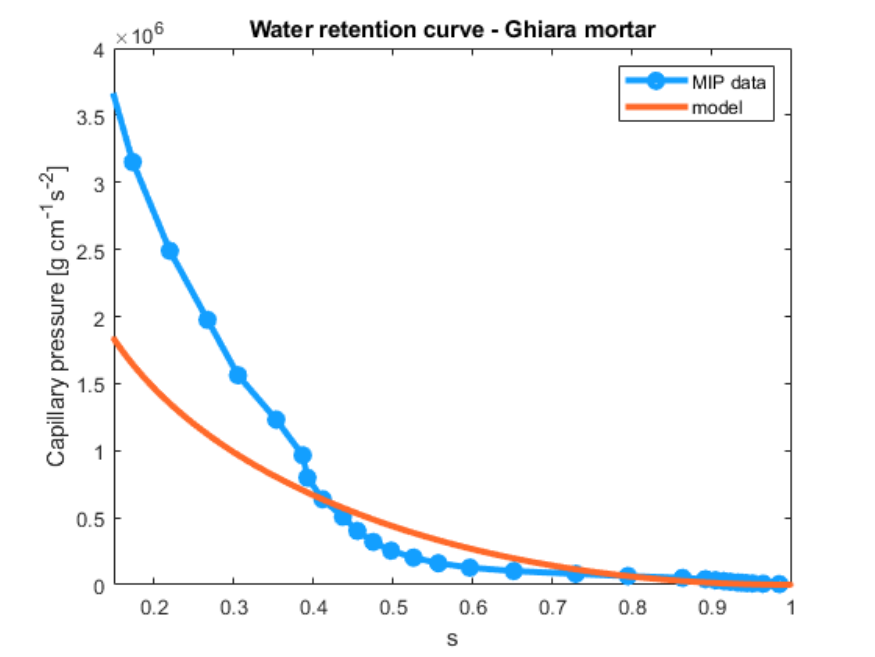}\quad \includegraphics[scale=0.4]{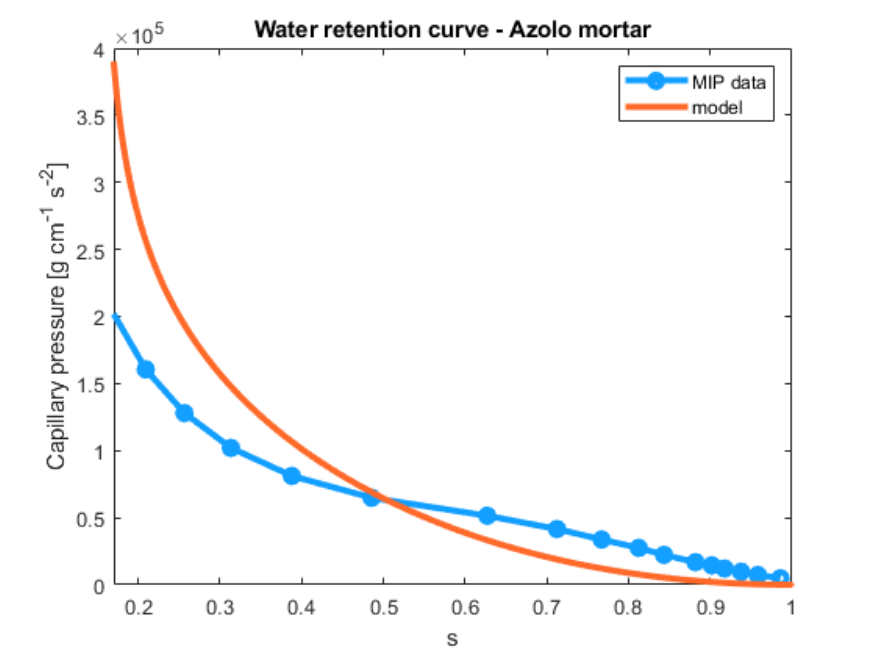}
    \caption{Left Panel: Plot of water retention curve for ghiara mortar. Experimental data (blue dotted line) VS curve profiles obtained numerically by the two models for capillary pressure $P_c(s)$, reported in \eqref{Pc} for $\alpha=0.25, c=$ 1.4e+06 (orange line). Right Panel: Plot of water retention curve for azolo mortar. Experimental data (blue dotted line) VS curve profiles obtained numerically by the two models for capillary pressure $P_c(s)$, reported in \eqref{Pc} for $\alpha=0.25, c=$1.98e+05 (orange line).} 
    \label{fig:MIP}
\end{figure}

Here we validate our procedure by comparing the model outcomes obtained for calibrated model parameters against MIP data described in paragraph \ref{sec:MIP}, by using the formulation of the capillary pressure function reported in \eqref{Pc}. In particular, in Fig. \ref{fig:MIP} are reported the profiles of the experimental  water retention curve and of the numerical one obtained with model \eqref{BkP} {using optimal parameters in Table \ref{table:3}.}

{It is worth noting that the superposition between the water retention curve derived from experimental data and that obtained by the mathematical model is used here as a qualitative and only partially quantitative comparison}. Indeed, our aim is to assess the order of magnitude of the capillary pressure values obtained from MIP experiment. In this way we can verify that model parameters obtained in the calibration procedure described in paragraph \ref{sec:calib} allow us to get values of $P_c(s)$ that are in accordance with experimental data, as shown for both materials in Fig. \ref{fig:MIP}.

%
\newpage
%


\begin{figure}
    \centering
    \includegraphics[scale=0.4]{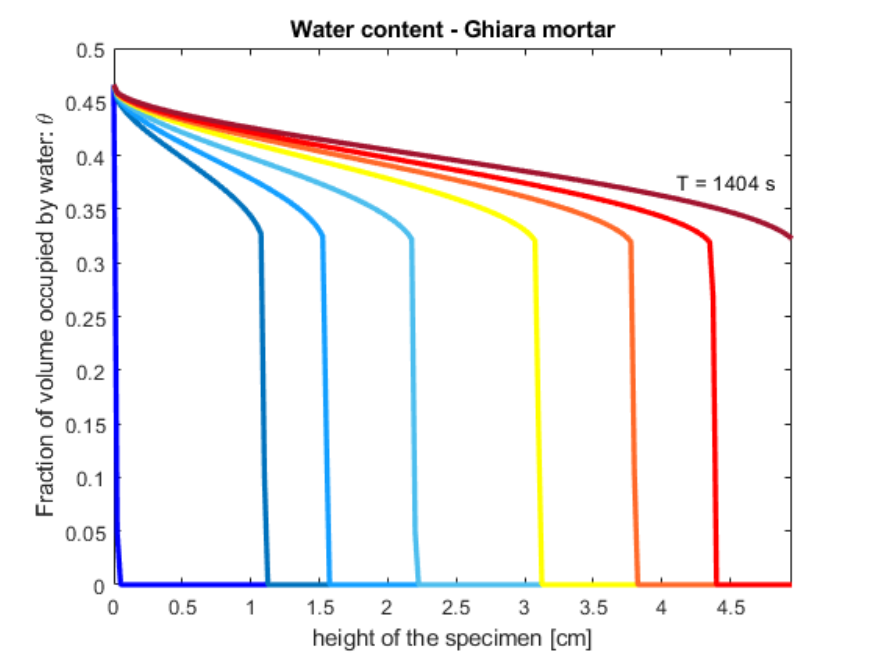}\quad \includegraphics[scale=0.4]{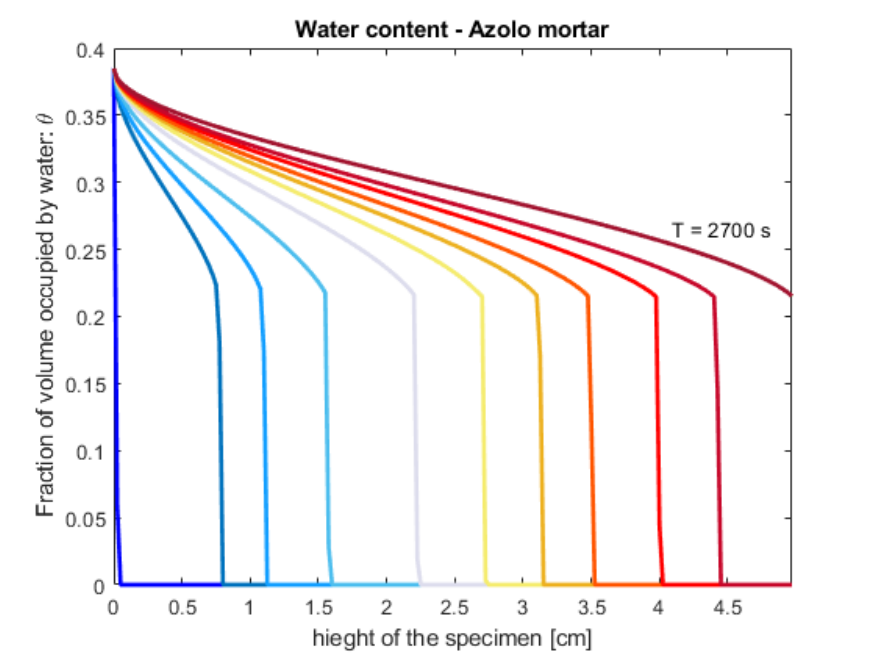}
    \caption{Left Panel: Plot of the water saturation as a function of length at different times during imbibition obtained numerically for ghiara mortar. Right Panel: Plot of the water saturation as a function of length at different times during imbibition obtained numerically for azolo mortar.}
    \label{fig:theta}
\end{figure}

\subsection{Sensitivity analysis}

In the present Section we perform a sensitivity analysis that explores the space of parameter values in which the mathematical problem is set, so that the effects of single parameters on the
dynamics can be evaluated.
{In the model \eqref{BkP} there are 6 parameter values that strongly affect
model outcomes in terms of the quantity of water absorbed by the specimen across time. Specifically, here we compute the error functional $E_2$ defined in \eqref{err_fun} and obtained by changing in turn only one parameter keeping the other fixed, i.e. one factor-at-time (OAT) procedure. Then, we observed that error functional $E_2$ is strictly related to the variation of parameter values and we verified that the optimal parameters obtained by the fitting procedure correspond to the minimum values of the error functionals, as expected.}

\begin{figure}
    \centering
\includegraphics[width=0.3\linewidth]{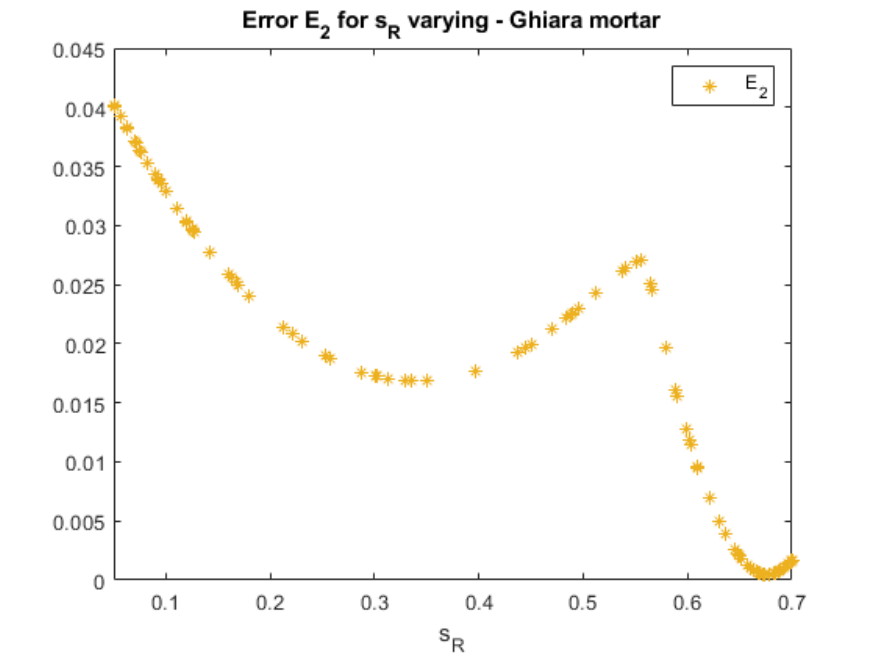}\quad\includegraphics[width=0.3\linewidth]{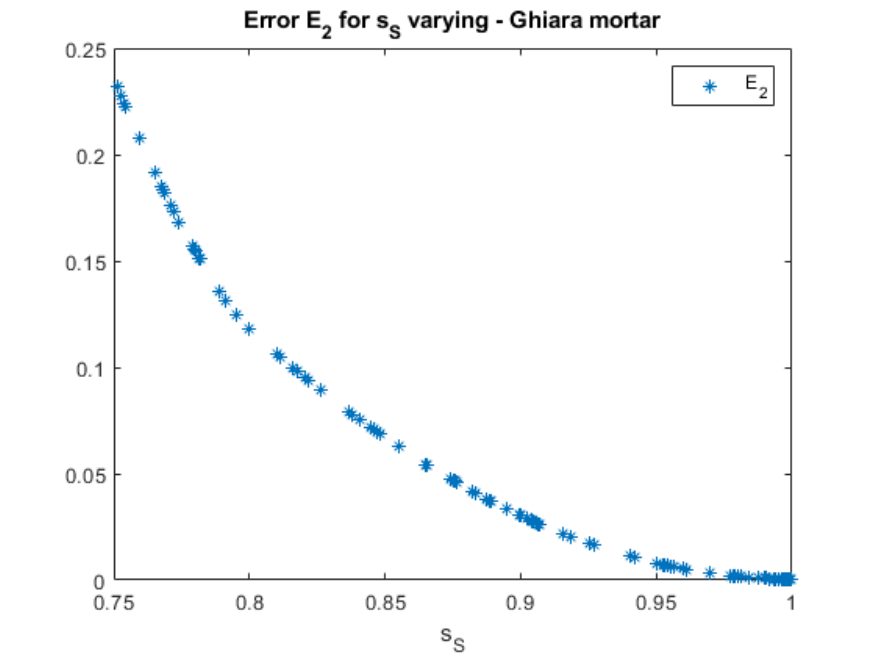}\quad\includegraphics[width=0.3\linewidth]{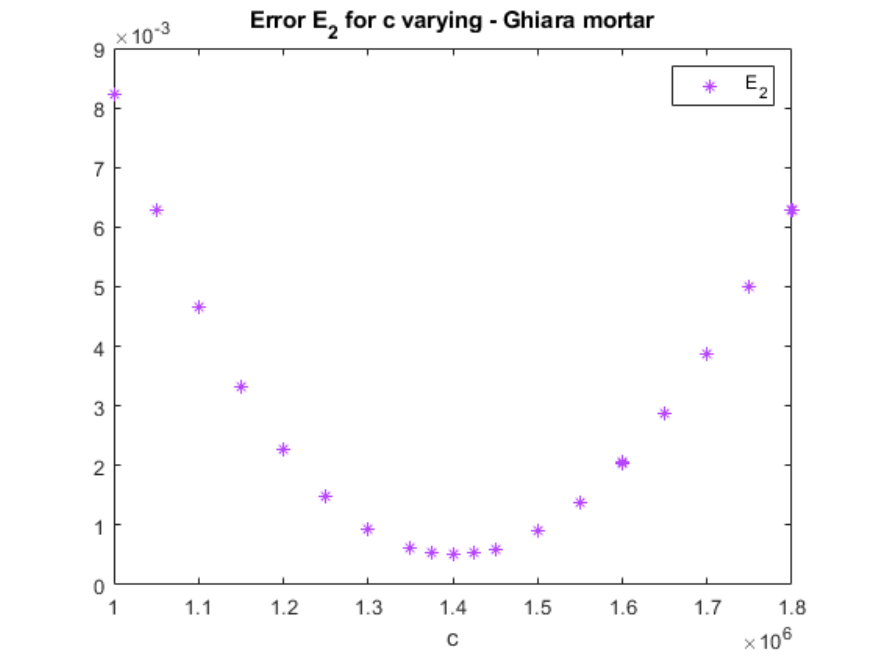}
\includegraphics[width=0.3\linewidth]{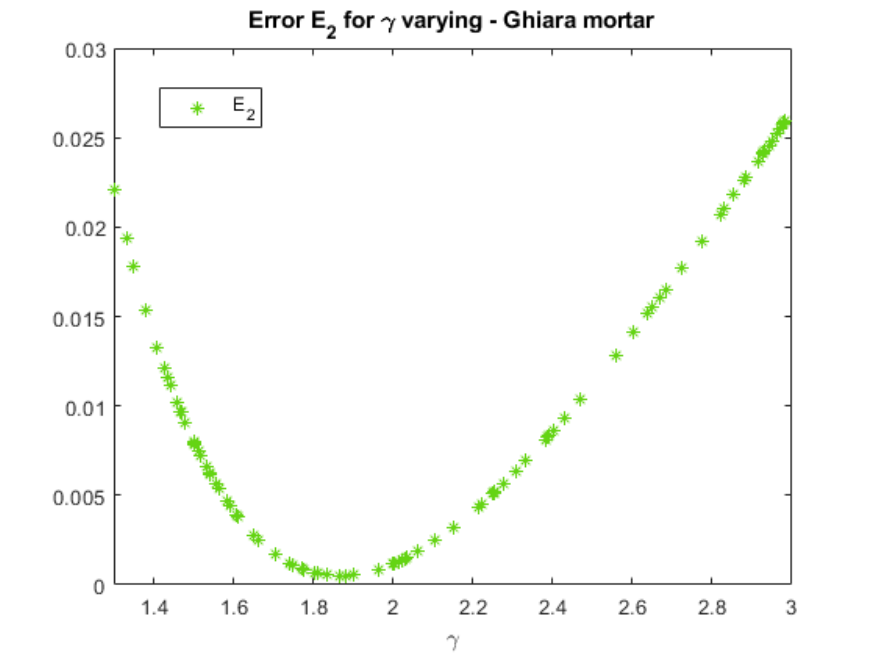}\quad\includegraphics[width=0.3\linewidth]{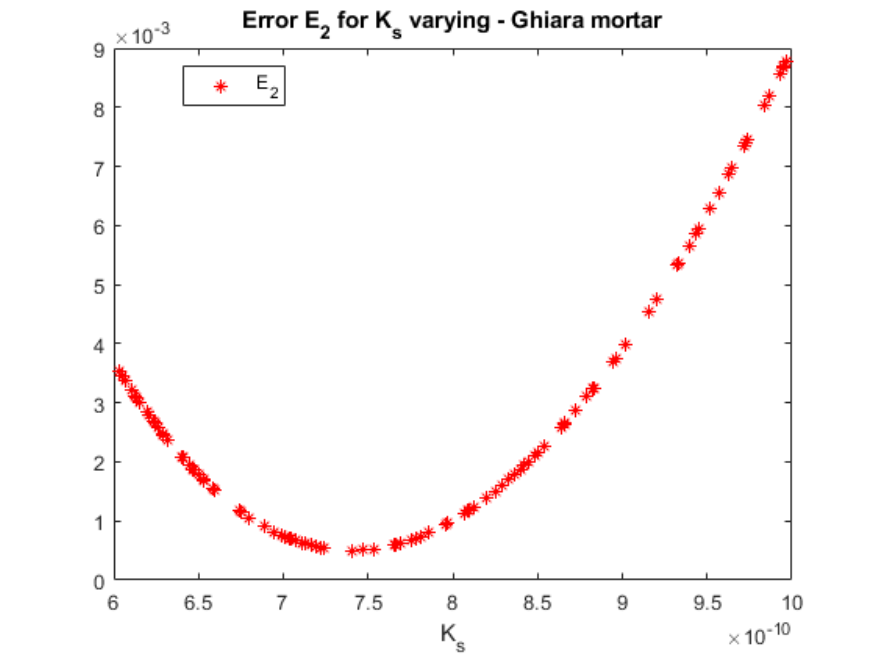}\quad\includegraphics[width=0.3\linewidth]{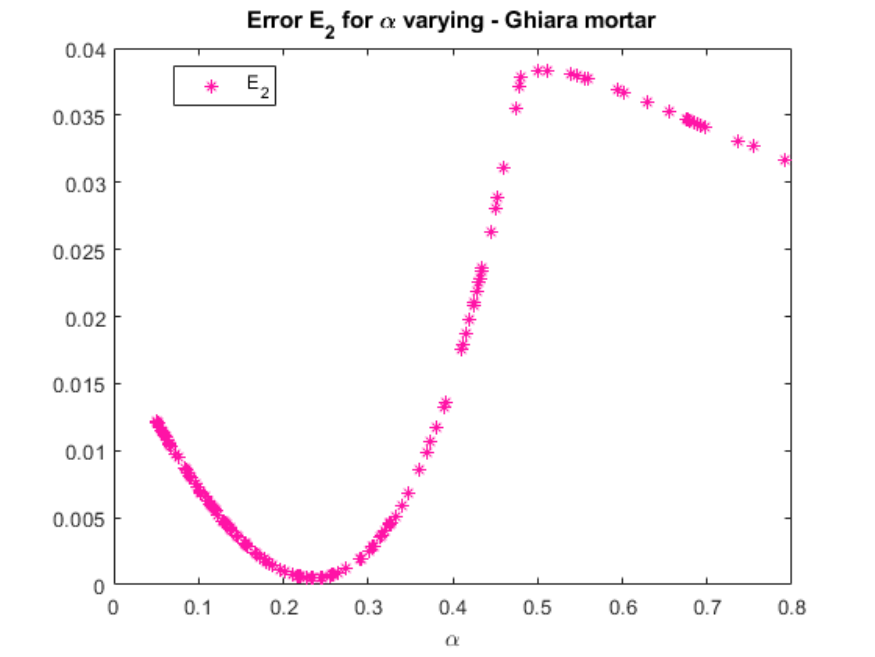}
    \caption{Local sensitivity analysis with respect to model parameters for ghiara mortar: plots of the error functional $E_2$.}
    \label{fig:errgh}
\end{figure}

\begin{figure}
    \centering
\includegraphics[width=0.3\linewidth]{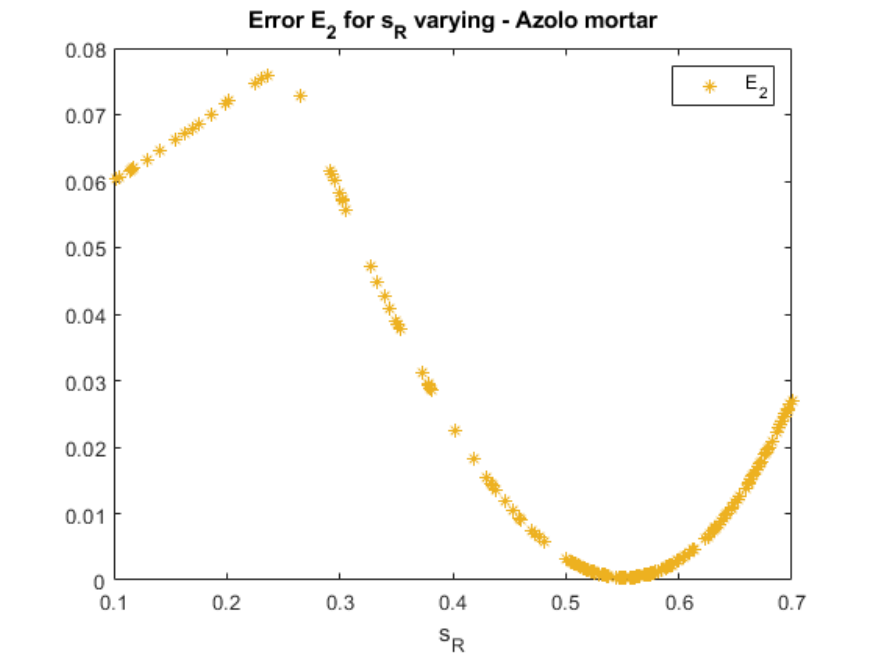}\quad\includegraphics[width=0.3\linewidth]{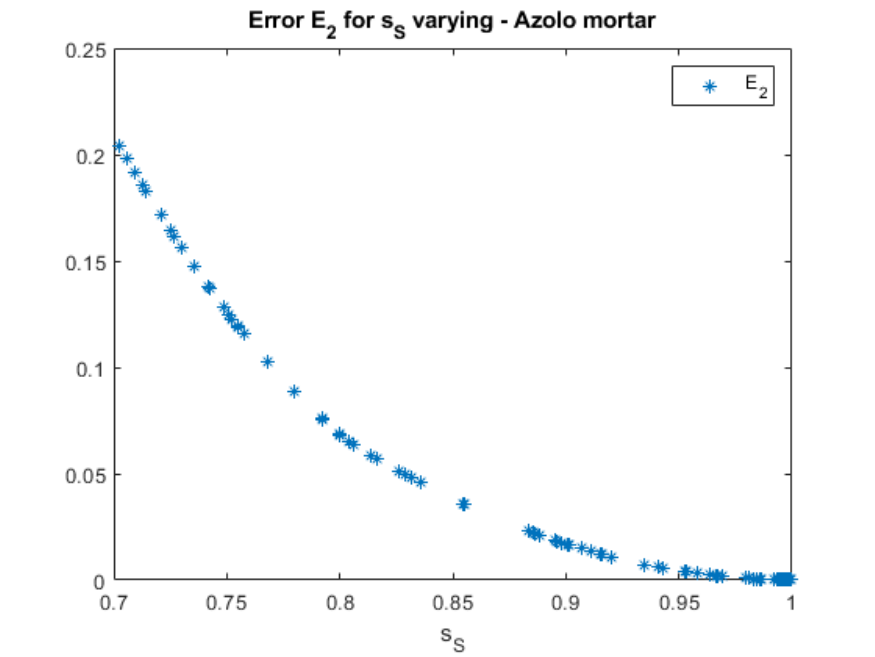}\quad\includegraphics[width=0.3\linewidth]{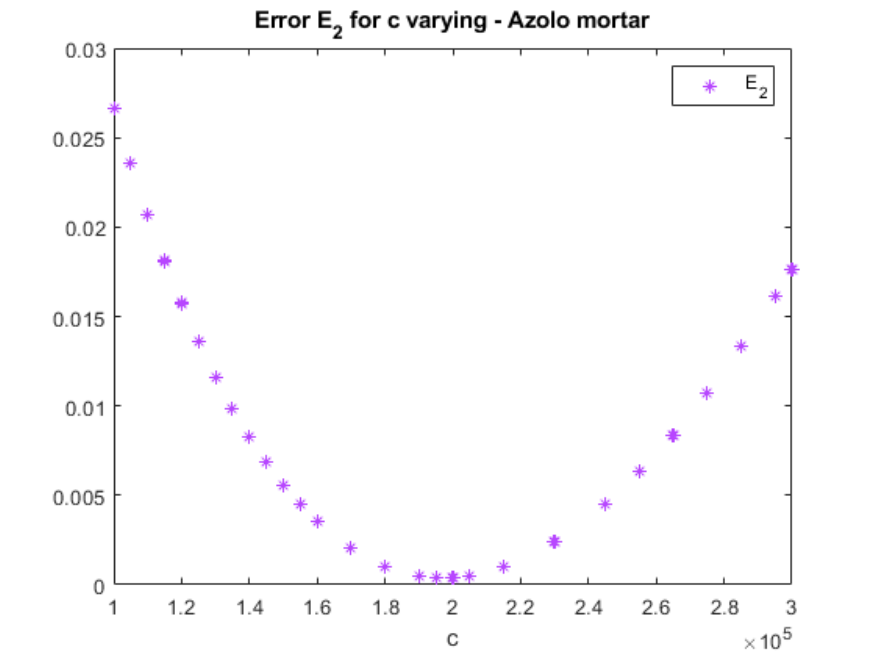}
\includegraphics[width=0.3\linewidth]{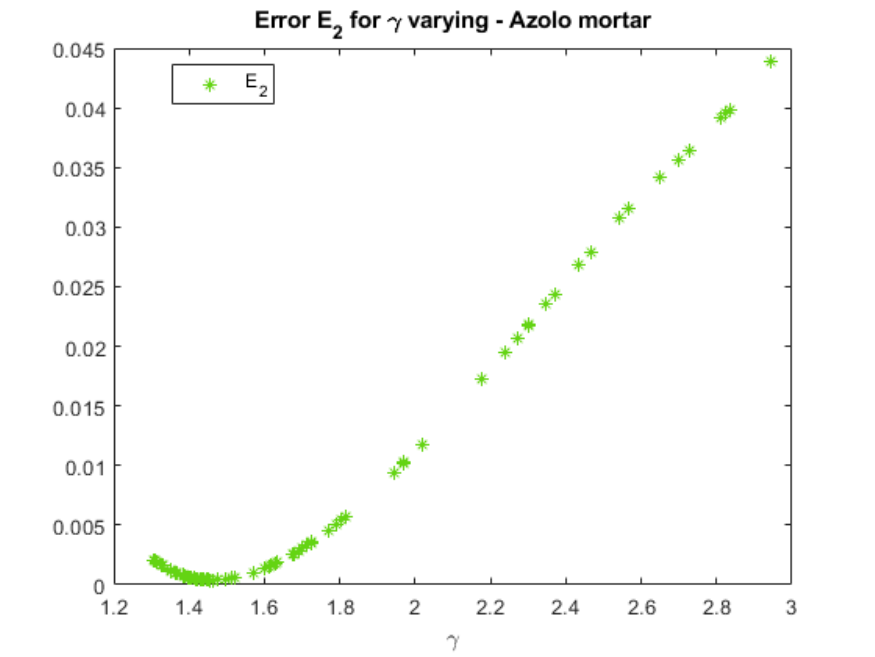}\quad\includegraphics[width=0.3\linewidth]{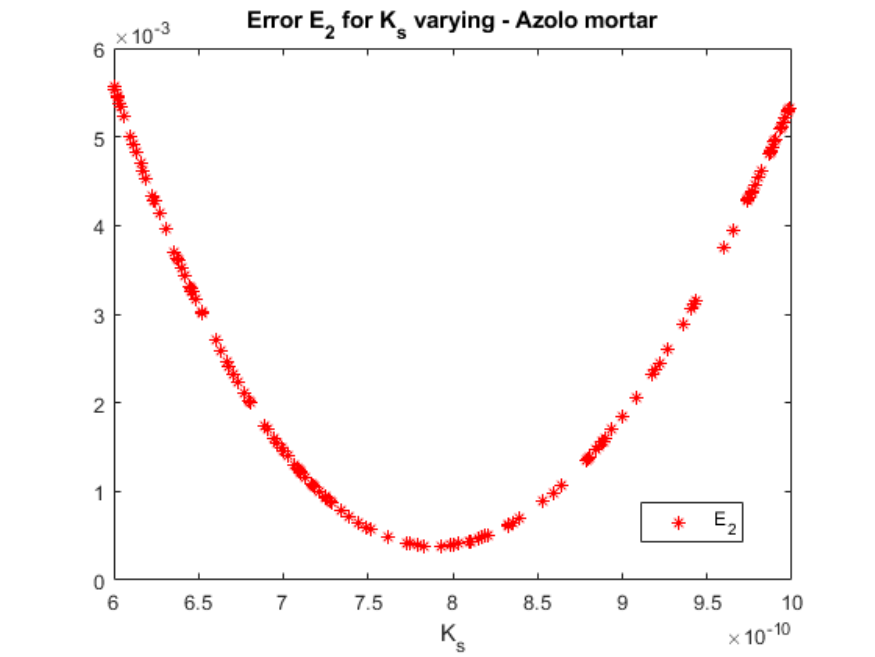}\quad\includegraphics[width=0.3\linewidth]{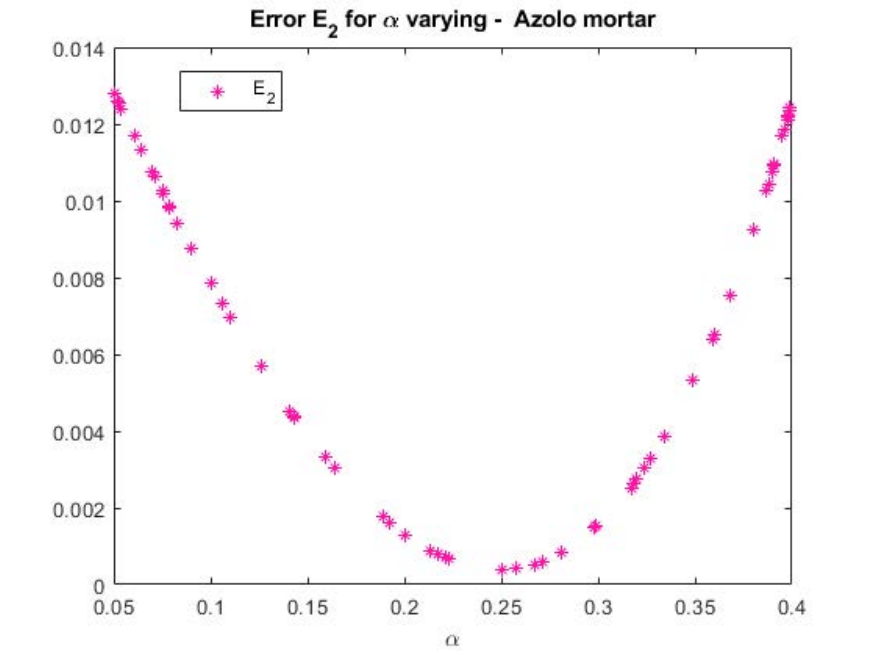}
    \caption{Local sensitivity analysis with respect to model parameters for azolo mortar: plots of the error functional $E_2$.}
    \label{fig:erraz}
\end{figure}

\section{Discussion of results and methodology}\label{sec:disc}

{From the results reported in Section \ref{sec:calib} and \ref{sec:valid}, we can see that water diffuses more rapidly in ghiara mortar respect to azolo one. Indeed, the diffusion coefficient $D_{kP}$ for ghiara mortar is about 4 times greater, as shown in Table \ref{table:3}. As a consequence, the water fills in the pores more rapidly in ghiara mortar, as can be seen from the profiles of water $\theta$ within the pores across time reported in Fig. \ref{fig:theta} along the specimen's heights. Then, the liquid in ghiara mortar reaches the sample's top more rapidly than in azolo one. 
Indeed, the water reaches the top of the specimen at time $T=2700 \ s$ (about 45 mins) for azolo mortar, while for ghiara mortar the top is reached at a shorter time $T= 1404 \ s$ (about 23 mins). This seems to be in accordance with the experimental data \cite{belfiore}: it can be seen that the imbibition curve reported in the left panel of Fig. \ref{fig:ghiara} in its initial part, i.e. until the first timepoint (10 mins), has a steepest slope with respect to the curve profile in the left panel of Fig. \ref{fig:azolo}.}

 Conversely, from 10 to 90 mins, the slope of imbibition curve of azolo mortar appears to be slightly steepest than that of ghiara mortar. This unlike behavior at different times can be justified by the fact that at the beginning of the test the ghiara mortar tends to absorb water more quickly due to its greater number of pores (its open porosity found by the MIP analysis is in fact higher, i.e. 46.6\% vs. 38.5\%, see Table \ref{table:phys}). However, it is also characterized by a more complex pore geometry with respect to azolo mortar and this, at a later time, causes in such a specimen a higher decrease in diffusive transport resulting from convolutions of the flow paths through porous media. Such assumption is supported by the "tortuosity" ($\tau$) values obtained from MIP analysis of the two mortar types analyzed in \cite{belfiore} and reported in Table \ref{table:phys}. More specifically, tortuosity is a parameter characterizing the geometry and length of interconnected phases. In porous media, it is defined as the ratio of the actual length of the flow-path divided by a straight line length in the direction of flow \cite{back}. The value $\tau$ of the two mortar specimens, as reported in Table \ref{table:phys}, is 7.6 for azolo and 9.9 for ghiara, which is in accordance with the above observations.
The final quantity of adsorbed water found by the mathematical model resulted to be higher in the ghiara mortar (2 $g/cm^2$) than in the azolo one (1.5 $g/cm^2$), even in this case fitting well the experimental data. Such a higher water content absorbed by the ghiara mortar can be again justified by its higher open porosity.

From the computational point of view, the numerical algorithm based on the simpler absorption model \eqref{NN} has two great advantages: on one hand, the running time is much lower respect to the numerical algorithm based on the model \eqref{BkP}; on the other hand the calibration procedure is more easily accomplished since there are only 3 parameters to be calibrated. {Both absorption models \eqref{NN} and \eqref{BkP} are able to reproduce quite well the phenomena observed experimentally, as can be seen from the values of errors obtained for optimal parameters and reported in Table \ref{table:err}. However, it can be noticed that the error values $E_2$ are less than $E_1$ for both materials, see Table \ref{table:err}. Moreover, it is worth to note that the absorption model \eqref{BkP} provides a deeper understanding on the permeability properties and capillary pressure of porous materials, since it provides the mathematical formulation of the permeability function as well as of the capillary pressure. For this reason, it allows us to have a validation of the model against experimental MIP data and permeability data, when available. }

In this case, the experimental data of permeability at saturation was not available, then the value of $K_s$ was obtained from the numerical calibration of the model by using data at our disposal. We point out that the value obtained for $K_s$ is of order of magnitude of $10^{-10} cm^2$ for both materials and thus falls within the ranges of values reported in literature for the material composing the mortars (i.e. lime, basalts and igneous rocks), see Fig. 1 in \cite{mappa_perm}. This result is qualitatively in accordance with the experimental findings in \cite{belfiore}, where the water vapor permeability test showed a similar behavior for the two mortars with only negligible differences between them.

\section{Conclusions and future perspectives}\label{concl}
In the present work we calibrated and validated a new mathematical model in order to describe the water absorption properties of two lime-based materials of the historic built heritage of Catania, i.e. ghiara and azolo mortars. The fitting procedure showed the effectiveness of our modelling approach, since we were able to simulate the amount of water absorbed by the specimen across time in accordance with imbibition data and to reproduce numerically the behavior of water retention curve derived from MIP data as well.
More in detail, we reproduced the imbibition curve with an average percentage error of about 0.01-0.02\%.

A validation of our modelling was possible by comparing the water retention curve obtained by the numerical algorithm against curves computed from MIP data.
In the future we aim at investigating further aspects both from the experimental and mathematical point of view.
From the experimental point of view we are interested in repeating the experiment of water imbibition for a longer observation time and applying protective treatments to the specimen made of ghiara and azolo mortars in order to see how the absorption properties change. 

From the mathematical viewpoint we will study the behavior of the models for water absorption here considered in comparison with classical models and we will apply them to different porous materials. Moreover, we will make a rigorous and systematic analysis of the effects of each model parameter on the numerical solution.
To this aim, in forthcoming works we will apply both theoretical methods and computational techniques of statistical inference to analyze the effect of the variation of model parameters on model outcomes.

\section{Acknowledgements.} G. B. is a member of the Gruppo Nazionale Calcolo
Scientifico-Istituto Nazionale di Alta Matematica (GNCS-INdAM). She is involved in PRIN project MATHPROCULT Prot. P20228HZWR, CUP J53D23015940001. G.B. is in the PNRR Project H2IOSC CUP B63C22000730005, financed by European Union - NextGenerationEU PNRR Mission 4, “Instruction and Research” - Component 2- Investment line 3.1. {G. B. and C. M. B. are in Project PE0000020 CHANGES - CUP: B53C22003890006 ,  NRP Mission 4 Component 2 Investment 1.3, Funded by the European Union - NextGenerationEU" under the Italian Ministry of University and Research (MUR).}

\end{document}